\documentclass{gtart_a}
\pdfoutput=1
\usepackage{graphicx}

\title{Siegel--Veech constants in $\mathcal{H}(2)$}

\author{Samuel Leli\`evre}
\givenname{Samuel}
\surname{Leli\`evre}
\address{Mathematics Institute\\
University of Warwick\\\newline
Coventry CV4 7AL\\UK}
\email{lelievre@maths.warwick.ac.uk}
\urladdr{http://carva.org/samuel.lelievre/}

\volumenumber{10}
\issuenumber{}
\publicationyear{2006}
\papernumber{27}
\lognumber{0588}
\startpage{1157}
\endpage{1172}
\MR{}
\Zbl{}

\arxivreference{math.GT/0503718}  


\received{30 March 2005}
\revised{23 February 2006}
\accepted{24 May 2006}
\published{12 September 2006}
\publishedonline{12 September 2006}
\proposed{Walter Neumann}
\seconded{Benson Farb, Tobias Colding}
\corresponding{}
\editor{}
\version{}

\AtBeginDocument{\let\bar\wbar\let\tilde\wtilde}
\def\S\,{Section }

\newcommand*{\newconcept}[1]{\textit{#1\/}}

\newcommand*\HH{\mathbf H}
\newcommand*\T{\mathbf T}

\newcommand*\cH{{\mathcal{H}}}
\newcommand*\cA{{\mathcal{A}}}
\newcommand*\cB{{\mathcal{B}}}
\newcommand*\cC{{\mathcal{C}}}
\newcommand*\cD{{\mathcal{D}}}
\newcommand*\cE{{\mathcal{E}}}
\newcommand*\cU{{\mathcal{U}}}

\DeclareMathOperator\SL{SL}
\DeclareMathOperator\PSL{PSL}
\DeclareMathOperator\vol{vol}
\DeclareMathOperator\lcm{lcm}
\DeclareMathOperator\cw{cw}
\newcommand*{\lquotient}[2]{\left.\raisebox{-0.2ex}{$#2$}%
  \backslash\raisebox{0.2ex}{$#1$}\right.}
\newcommand*{\rquotient}[2]{\left.\raisebox{0.2ex}{$#1$}%
  /\raisebox{-0.2ex}{$#2$}\right.}
\newcommand\mat[4]{
  \bigl( \begin{smallmatrix} #1&#3\\ #2&#4 
  \end{smallmatrix}\bigr)
  }
\newcommand\hp[1]{
  \mat{1}{0}{#1}{1}%
  }
\newtheorem{Theorem}{Theorem}
\newtheorem{Lemma}{Lemma}
\makeautorefname{Lemma}{Lemma}
\makeautorefname{Theorem}{Theorem}

\theoremstyle{remark}
\newtheorem*{Remark}{Remark}

\begin{document}

\begin{asciiabstract}
Abelian differentials on Riemann surfaces can be seen 
as translation surfaces, which are flat surfaces with 
cone-type singularities.
Closed geodesics for the associated flat metrics form 
cylinders whose number under a given maximal length 
generically has quadratic asymptotics in this length,
with a common
coefficient constant for the quadratic asymptotics called a 
Siegel--Veech constant which is shared by almost all
surfaces in each moduli space of translation surfaces.

Square-tiled surfaces are specific translation surfaces which
have their own quadratic asymptotics for the number of cylinders of
closed geodesics. It is an interesting question whether, as n tends to infinity,
the Siegel--Veech constants of square-tiled surfaces with n tiles
tend to the generic constants of the ambient moduli space.
We prove that this is the case in the moduli space H(2) of
translation surfaces of genus two with one singularity.
\end{asciiabstract}

\begin{htmlabstract}
<p class="noindent">
Abelian differentials on Riemann surfaces can be seen as translation
surfaces, which are flat surfaces with cone-type singularities.
Closed geodesics for the associated flat metrics form cylinders whose
number under a given maximal length was proved by Eskin and Masur to
generically have quadratic asymptotics in this length, with a common
coefficient constant for the quadratic asymptotics called a
Siegel&ndash;Veech constant which is shared by almost all
surfaces in each moduli space of translation surfaces.
</p>
<p class="noindent">
Square-tiled surfaces are specific translation surfaces which
have their own quadratic asymptotics for the number of cylinders of
closed geodesics. It is an interesting question whether the Siegel&ndash;Veech constant of a
given moduli space can be recovered as a limit of individual constants
of square-tiled surfaces in this moduli space.
We prove that this is the case in the moduli space H(2) of
translation surfaces of genus two with one singularity.
</p>
\end{htmlabstract}

\begin{abstract}
Abelian differentials on Riemann surfaces can be seen as translation
surfaces, which are flat surfaces with cone-type singularities.
Closed geodesics for the associated flat metrics form cylinders whose
number under a given maximal length was proved by Eskin and Masur to
generically have quadratic asymptotics in this length, with a common
coefficient constant for the quadratic asymptotics called a 
Siegel--Veech constant which is shared by almost all
surfaces in each moduli space of translation surfaces.

Square-tiled surfaces are specific translation surfaces which
have their own quadratic asymptotics for the number of cylinders of
closed geodesics. It is an interesting question whether the Siegel--Veech constant of a
given moduli space can be recovered as a limit of individual constants
of square-tiled surfaces in this moduli space.
We prove that this is the case in the moduli space $\mathcal{H}(2)$ of
translation surfaces of genus two with one singularity.
\end{abstract}

\maketitle

\section{Introduction}

\subsection{Geodesics on the torus}

On the standard torus $\T^2 = \rquotient{\R^2}{\Z^2}$, the number
$N(L)$ of maximal families of parallel simple closed geodesics of
length not exceeding $L$ is well-known (and easily seen) to grow
quadratically in $L$, with
$$
N(L)\sim\frac{1}{2\,\zeta(2)} \cdot \pi L^2,
$$
which is half of the asymptotic for the number of primitive
lattice points in a disc of radius $L$.
The factor one-half comes from counting unoriented rather than
oriented geodesics.

By convention, the corresponding \emph{Siegel--Veech constant} is
$\unpfrac{1}{2\,\zeta(2)}$.
Note that it is the coefficient of $\pi L^2$,
not of $L^2$, in the asymptotic.

Marking the origin of the torus (
artificially considering it as
a singularity or saddle), the number of geodesic segments of length at most $L$ joining the
saddle to itself coincides with the number of
families of simple closed geodesics.

\subsection{Geodesics on translation surfaces}

It is a standard fact that abelian differentials on Riemann
surfaces can be seen as translation surfaces.
On translation surfaces of genus at least $2$, countings of closed
geodesics or saddle connections similar to those just described for
the torus can be made.

There, the countings of saddle connections and of cylinders of simple
closed geodesics do not coincide, but their growth rates remain
quadratic. This is made more precise by several related results.

Masur \cite{Ma88,Ma90} proved that for every translation surface,
there exist positive constants $c$ and $C$ such that the counting
functions of saddle connections and of maximal cylinders of closed
geodesics satisfy
$$
c \cdot \pi L^2 \leqslant N_{\mathrm{cyl}}(L) \leqslant N_{\mathrm{sc}}(L)
\leqslant C \cdot \pi L^2
$$
for large enough $L$.

Veech \cite{Ve89} proved that on a square-tiled surface (and on any
Veech surface) there are in fact \emph{exact quadratic asymptotics};
Gutkin and Judge \cite{GuJu} gave a different proof.

Another proof for the upper quadratic bounds for $N_{\mathrm{cyl}}(L)$
and $N_{\mathrm{sc}}(L)$ was given by Vorobets \cite{Vo97}.
Eskin and Masur \cite{EM} gave yet another one and proved that for
each connected component of each stratum of each moduli space of
normalised (area $1$) abelian or quadratic differentials, there
are constants $c_{\mathrm{sc}}$ and $c_{\mathrm{cyl}}$ such that
\emph{almost every surface} in the component has $N_{\mathrm{sc}}(L)
\sim c_{\mathrm{sc}} \pi L^2$ and $N_{\mathrm{cyl}}(L) \sim c_{\mathrm{cyl}}
\pi L^2$.

It is an interesting open problem whether \emph{all} translation
surfaces have exact quadratic asymptotics for countings of saddle
connections and of cylinders of closed geodesics.

The particular constants for many Veech surfaces have been computed
explicitly by Veech \cite{Ve89}, Vorobets \cite{Vo97}, Gutkin and Judge
\cite{GuJu}, and Schmoll \cite{Schmo}.  Constants for some families of
non-Veech surfaces were also given by Eskin, Masur and Schmoll \cite{EMS}
and Eskin, Marklof and Morris \cite{EMWM}.
The generic constants for the connected components of all strata of
abelian differentials were computed by Eskin, Masur and Zorich
\cite{EMZ}.

In general, the particular constants for Veech surfaces do not
coincide with the generic constants of the strata where they live.

There is another subtle difference between Veech surfaces and
generic surfaces.
Define cylinders as \newconcept{regular} if their boundary
components both consist of a single saddle connection.
In any connected component of stratum in genus at least $2$, a
generic surface has no irregular cylinders while on Veech surfaces,
countings of irregular cylinders have quadratic asymptotics.

However we will prove that on the stratum
$\cH(2)$ of translation surfaces of genus $2$ with one singularity,
the individual `quadratic constants' for \emph{regular} cylinders on square-tiled surfaces 
retreive the generic Siegel--Veech constant of $\cH(2)$ as a
limit.
See \fullref{thm:cv:svc:reg:cyl} in \fullref{sec:stratum:h2} for a
precise statement.

\subsection{Setting and main result}
\label{sec:stratum:h2}

In this paper, we are concerned with the stratum $\cH(2)$ consisting
of genus--$2$ abelian differentials with a double zero, or in other
words, translation surfaces of genus $2$ with one singularity (of angle
$6\pi$).
\begin{Theorem}
\label{thm:cv:svc:reg:cyl}
Consider a sequence $S_n$ of area--$1$ surfaces in $\cH(2)$, each tiled
by some prime number $p_n$ of squares, with $p_n\to \infty$.
Then the constants in the quadratic asymptotics for regular cylinders of
closed geodesics on the surfaces $S_n$ tend to
$\punfrac{10}{3}(\unpfrac{1}{2\zeta(2)})$, the Siegel--Veech constant of
$\cH(2)$ for cylinders of closed geodesics.
\end{Theorem}

\begin{Remark}
It is possible to adapt our calculations to show that the constants in
the quadratic asymptotics for \emph{irregular} cylinders of closed
geodesics on the surfaces $S_n$ in the theorem tend to $0$, so that
the constants in the quadratic asymptotics for \emph{all} cylinders
(both regular and irregular) tend to the generic constant of the
stratum $\cH(2)$ as well.
\end{Remark}

\begin{Remark}
We believe that the assumption that the number of squares tiling the
surfaces is prime is unnecessary, but we have not yet been able to
adapt the calculations to show the convergence of Siegel--Veech
constants in the case of nonprime numbers of tiles.
\end{Remark}

The proof of the theorem relies on fine estimates presented in
\fullref{sec:simpler:sum}.

Pierre Arnoux pointed out to us the analogy to a result of C\,Faivre
on L\'evy constants of quadratic numbers.  See Faivre \cite{F} or Dal'Bo and Peign\'e \cite{DP}.

\subsection{Acknowledgments}

The author wishes to thank Anton Zorich for guiding him to this
problem, Pascal Hubert, Jo\"el Rivat and Emmanuel Royer for useful
conversations, and C\'ecile Dartyge and G\'erald Tenenbaum who helped
him with the estimates in \fullref{sec:simpler:sum}. He also thanks
the referee for useful remarks.

This research was carried out in Montpellier for the most part; some
intuition was gained from computer calculations (programmed in Caml
Light) run using the Medicis server at \'Ecole polytechnique.

\section{Preliminaries}

The stratum $\cH(2)$ is the simplest stratum of abelian differentials
after the (well-understood) stratum of abelian differentials on tori.
As every stratum, it admits a natural $\SL(2,\R)$ action, and we will
recall here some facts concerning the orbits of certain special points
of $\cH(2)$, the square-tiled surfaces.

A square-tiled surface is a ramified translation cover of the standard
torus with only one branch point.  The number of square tiles is the
number of sheets of the covering or the degree of the corresponding
covering map to the standard torus.  A square-tiled surface is called
\emph{primitive} if this covering map does not factor through a
covering of a larger torus with only one branch point.

\subsection{Orbits of square-tiled surfaces}
\label{sec:orbits}

By a theorem of McMullen \cite{McDS}, in $\cH(2)$, 
primitive $n$--square-tiled surfaces for $n > 3$ are in a single $\SL(2,\R)$--orbit
if $n$ is even, and in exactly two $\SL(2,\R)$--orbits if $n$ is odd (see
Hubert and Leli\`evre \cite{HL1} for the prime $n$ case).
We will denote these orbits by $\cA_n$ and $\cB_n$ for odd $n$ and by
$\cE_n$ for even $n$.

The integer points in these orbits are primitive $n$--square-tiled
surfaces, and they form $\SL(2,\Z)$--orbits which we will denote
respectively by $A_n$, $B_n$ and $E_n$.
The number of primitive $n$--square-tiled surfaces in $\cH(2)$ is thus
the cardinality of $E_n$ when $n$ is even and the sum of the
cardinalities of $A_n$ and $B_n$ when $n$ is odd.  This number is
given in \cite[Lemma~4.11]{EMS} to be asymptotic to
$$
\frac{3}{8}\,n^3\,\prod_{p \mid n}\Big(1 - \frac{1}{p^2}\Big).
$$
Formulas for the separate countings of $A_n$ and $B_n$ conjectured
by Hubert and Leli\`evre \cite{HL1} are established by Leli\`evre and Royer \cite{LeRo} to be
respectively
$$
a_n = \frac{3}{16}\,(n-1)\,n^2\,\prod_{p \mid n}\Big(1 - \frac{1}{p^2}\Big)
\quad \text{and} \quad
b_n = \frac{3}{16}\,(n-1)\,n^2\,\prod_{p \mid n}\Big(1 - \frac{1}{p^2}\Big).
$$
If $n$ tends to infinity within the set of prime numbers, both $a_n$ 
and $b_n$ are asymptotic to $\punfrac{3}{16}\,n^3$.

The natural definition of primitive $n$--square-tiled surfaces gives
them area $n$ (each square tile has area $1$), but it is sometimes
useful to consider the corresponding unit area surfaces by applying
the natural projection from $\cH(2)$ to the unit hyperboloid $\cH_1(2)$.

\subsection{Cusps}
\label{sec:cusps}

Each square-tiled surface in the stratum $\cH(2)$ decomposes into either
one or two horizontal cylinders, and can be given as coordinates the
heights, widths and twist parameters of these cylinders; see \cite{EMS}
or \cite{HL1}.
Here we are interested in \emph{regular} cylinders of closed
geodesics, which exist only in \emph{two-cylinder} decompositions
(in one-cylinder decompositions, the unique cylinder has three saddle 
connections on each boundary component).

The decompositions into cylinders provide a way to parametrise
square-tiled surfaces (by the heights, widths and twist parameters of 
their cylinders). 
These parameters are very convenient to describe the action of $\cU =
\{\mat{1}{0}{n}{1}\colon n \in \Z\}$, which only affects the twist
parameters.

The following lemma puts together Lemmas 2.4, 2.5 and 3.1 of 
\cite{HL1}.  The notations $\wedge$ and $\vee$ are used for $\gcd$ and $\lcm$, 
respectively.

\begin{Lemma}
\label{cusp:lemma}
Let $S$ be a primitive $n$--square-tiled surface, and denote by $\cD$,
resp.\ D, its orbit under $\SL(2,\R)$, resp.\ $\SL(2,\Z)$.  Then $D$
is the set of primitive $n$--square-tiled surfaces in $\cD$ and the
cusps of $\cD$ are in bijection with the $\cU$--orbits in $D$.

If $S$ has two cylinders, with $h_i$, $w_i$ and $t_i$ ($i = 1,\ 2$) as
height, width and twist parameters, then its cusp width (the
cardinality of its $\cU$--orbit) is

$$\displaystyle \mathrm{cw}(S) =
\frac{w_1}{w_1 \wedge h_1} \vee \frac{w_2}{w_2 \wedge h_2}
\quad \Big(=\frac{w_1}{w_1 \wedge h_1} \times \frac{w_2}{w_2 \wedge h_2}
\text{ for prime } n\Big).
$$

The surface $S'$ with $h'_i = h_i$, $w'_i = w_i$, and $t'_i = t_i
\bmod (w_i \wedge h_i)$ is a ``canonical'' representative of the
$\cU$--orbit of $S$.
Each cusp thus has a unique representative with $0 \leqslant t'_i <
w_i \wedge h_i$.
\end{Lemma}

We also recall that each direction
of rational slope on a given square-tiled surface $S$ gives rise to a decomposition of $S$ in
cylinders of closed geodesics, and this direction can be
associated to one of the cusps of the $\SL(2,\R)$--orbit of $S$.

Note that these cusps can also be understood as cusps of $\Gamma(S)$,
the Veech group, or stabiliser under $\SL(2,\R)$ of $S$. 
Algebraically this means conjugacy classes of maximal parabolic 
subgroups; geometrically the `cusps' of the quotient surface 
$\lquotient{\HH}{\Gamma(S)}$.

\subsection{A formula for the constants}
\label{sec:formula}

Here, we establish a formula for the constants for which we will
compute estimates in \fullref{sec:asymptotics}.

\begin{Lemma}
\label{lemma:asymptotics:nreg}
The number $N_{\mathrm{reg}}(L)$ of regular cylinders of closed
geodesics of length at most $L$ on a unit area square-tiled surface
$S$ has the following asymptotics:
$$
N_{\mathrm{reg}}(L) \sim \frac{n}{\#D} \sum_{\substack{\cC_j\text{
two-cyl} \\ \text{cusp of }S}} \frac{\cw(\cC_j)}{w_1^2}
\,\frac{1}{2\zeta(2)}\,\pi L^2.
$$
\end{Lemma}

\begin{Remark}
Following tradition we write the asymptotic as a multiple of $\pi
L^2$ rather than just $L^2$ and write $(\unpfrac{1}{2\zeta(2)})\pi L^2$
instead of $\punfrac{3}{\pi} L^2$ to bring out the analogy with the
corresponding formula for the torus.
\end{Remark}

\begin{proof}
We deduce this formula from the material reviewed in
\fullref{sec:orbits}--\fullref{sec:cusps}, and from \S\,3 of
\cite{Ve89} to which we refer freely here both for notation and
results.

For any finite covolume subgroup $\Gamma$ of
$\PSL(2,\R)$, Veech introduces a complete set $\{\Lambda_j\}_{1\leqslant j \leqslant
r}$ of representatives of the maximal parabolic subgroups of $\Gamma$.
We will also refer to the cusps $\cC_j$.
He defines $\Lambda_0 = \{\hp{k}\colon k \in \Z\}$ which we denoted by
$\cU$.

Then for each $j$ he selects $\beta_j\in\SL(2,\R)$ conjugating
$\Lambda_0$ to $\Lambda_j$, ie, $\beta_j^{-1}\Lambda_j\beta_j =
\Lambda_0$.
When $\Gamma$ is the Veech group of a translation surface $S$, this
amounts to representing the cusp $j$ by the surface $\beta_j^{-1} S$.
Indeed, $\beta_j S$ has Veech group $\beta_j^{-1} \Gamma \beta_j$.

These `representatives' $\beta_j^{-1}S$ of the cusps have width $1$.
For square-tiled surfaces, $\Gamma(S)$ is always a subgroup of
$\SL(2,\Z)$, so it is also usual to conjugate inside $\SL(2,\Z)$ to
the group generated by some $\hp{k}$ rather than to $\Lambda_0$
itself, thus keeping track of the cusp width (the adequate $k$).

Let us illustrate the difference on an example.

Consider the surface $S$ pictured on the left of
\fullref{fig:surface} made of seven squares $s_1,\ldots, s_7$
forming a horizontal cylinder where the right edge of each $s_i$ is
glued to the left edge of $s_{i+1}$ (indices being understood modulo
$7$), and where the top edges of squares $s_1$ to $s_7$ are
respectively glued to the bottom edges of squares $s_3$ to $s_6$,
$s_1$ to $s_2$, and $s_7$.

Consider the direction of the first diagonal.  In this direction the
surface $S$ decomposes into cylinders of closed geodesics as
illustrated on the right of \fullref{fig:surface} (parallel sides
of same length are identified).

\begin{figure}[ht!]
\centering
\includegraphics{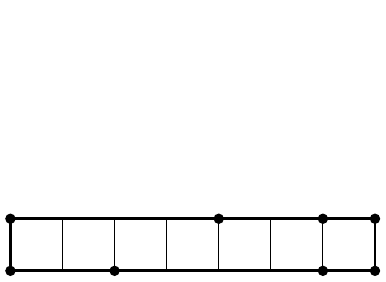}%
\qquad
\includegraphics{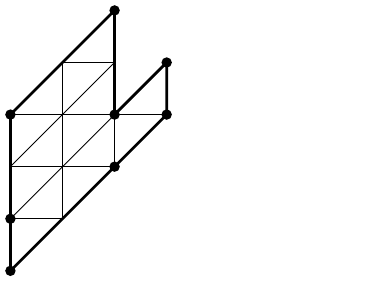}
\label{fig:surface}
\caption{}
\end{figure}

We get to the standard square-tiled representative of the cusp
corresponding to that direction by applying the matrix $M =
\mat{\ 1}{-1}{0}{1}$ (the matrix in $\SL(2,\Z)$ which sends $(1,1)$ to
$(1,0)$ and $(0,1)$ to itself).  A choice of $\beta_j^{-1}$ is
$\mat{1/\sqrt{3}}{0}{0}{\sqrt{3}}\cdot M$.  This sends $3\,{\cdot}\,(1,1)$ 
to $(1,0)$.

\fullref{fig:cusp:representative} represents
$\beta_j^{-1}\,S$ on the left and $M\,S$ on the right.

\begin{figure}[ht!]
\centering
\includegraphics{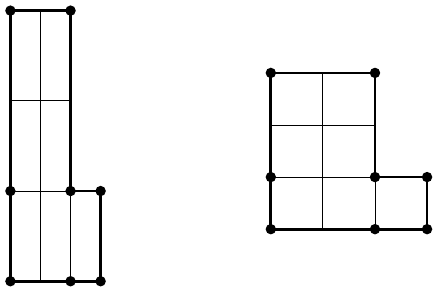}
\label{fig:cusp:representative}
\caption{}
\end{figure}

Veech defines $\xi_j$ to be the vector $\beta_j\,{\cdot}\,(1,0)$.  And for
each cylinder of closed geodesics in the direction of $\Lambda_j$,
calling $v$ the holonomy vector for this cylinder, he associates to the
cylinder the constant $c_j(v) =
\unfrac{\left\|\xi_j\right\|}{\left\|v\right\|}$.

In our notation, for a surface tiled by unit squares, if $v$ is the
holonomy vector of cylinder $i$ of cusp $\cC_j$, we have $c_j(v) =
\sqrt{\cw(\cC_j)}/w_i$, where $w_i$ is the width of this cylinder and
$\cw(\cC_j)$ the width of this cusp.

Veech's formula for the asymptotics \cite[formula (3.11)]{Ve89} 
is:
$$
N(L) \sim \vol(\lquotient{\HH}{\Gamma(S)})^{-1}
\Biggl(\sum_{j=1}^r\biggl(\sum_v c_j(v)^2\biggr)\Biggr) L^2.
$$

So the contribution of a given cylinder of a cusp $\cC_j$ to the
coefficient of the quadratic asymptotics is
$\vol(\lquotient{\HH}{\Gamma(S)})^{-1} c_j(v)^2$, where $v$ is the 
holonomy vector of this cylinder.

If we are concerned with regular cylinders of closed geodesics for 
square-tiled surfaces in $\cH(2)$, we need only consider cylinder $1$ 
of two-cylinder cusps.

The volume of the quotient $\lquotient{\HH}{\Gamma(S)}$ equals the
index of $\Gamma(S)$ in $\SL(2,\Z)$ times the volume of
$\lquotient{\HH}{\SL(2,\Z)}$; and the index equals the cardinality of
the $\SL(2,\Z)$--orbit $D$ of $S$, while the volume of
$\lquotient{\HH}{\SL(2,\Z)}$ is $\pi/3$.

The last thing to observe is the effect of scaling a surface.
Consider a surface $S$ where quadratic asymptotics $N(L) = c\cdot \pi
L^2$ hold, and scale $S$ by a scale factor $r$.  On $r S$, the
asymptotics become $N(L) = c \cdot \pi (L/r)^2$.

A square-tiled surface $S$ of area $1$ is a scaled-down version by
$\unfrac{1}{\sqrt{n}}$ of a surface tiled by $n$ unit squares; this scaling
changes the asymptotic by a factor $n$.

This completes the proof of the formula in \fullref{lemma:asymptotics:nreg}.
\end{proof}

Let us denote by $\tilde{c}(S)$ the quantity 
$$
\frac{n}{\#D} \sum_{\substack{\cC_j\text{ two-cyl} \\ \text{cusp of
}S}} \frac{\cw(\cC_j)}{w_1^2}.
$$
Our aim is now to prove that $\tilde{c}(S)$ tends to
$\unfrac{10}{3}$ as the number of square tiles of $S$ tends to infinity
staying prime. This will establish \fullref{thm:cv:svc:reg:cyl}.

As a first step for this, using the description of the two-cylinders
cusps in the orbits of square-tiled surfaces (see \cite{HL1}), and
renaming $w_1$, $w_2$, $h_1$, $h_2$ as $a$, $b$, $h$, $y$
respectively, we get:
\begin{itemize}
\item
for $S$ in orbit $A_n$, 
$$
\tilde{c}(S) = \frac{n}{a_n} \Bigl(\sum_{\substack{a,b,h,y
\geqslant 1 \\ a h + b y = n \\ h \wedge y = 1 \\ \text{$h$, $y$ odd}
\\ a < b}} \frac{a b}{a^2} + \frac{1}{2} \sum_{\substack{a,b,h,y
\geqslant 1 \\ a h + b y = n \\ h \wedge y = 1 \\ a \not \equiv b\!\!\mod
2 \\ h \not \equiv y\!\!\mod 2 \\ a < b}} \frac{a b}{a^2}\Bigr)
$$
\item
for $S$ in orbit $B_n$, 
$$
\tilde{c}(S) = \frac{n}{b_n}
\Bigl(\sum_{\substack{a,b,h,y\geqslant 1 \\ a h + b y = n \\ h 
\wedge y = 1 \\ \text{$a$, $b$ odd} \\ 
a < b}} \frac{a b}{a^2} + \frac{1}{2} \sum_{\substack{a,b,h,y
\geqslant 1 \\ a h + b y = n \\ h \wedge y = 1 \\ a \not \equiv b \!\mod
2 \\ h \not \equiv y \!\mod 2 \\ k < \ell}} \frac{a b}{a^2}\Bigr)
$$
\end{itemize}

The idea is to group two-cylinder cusps sharing the same parameters 
$w_1$, $w_2$, $h_1$, $h_2$.
Then the sum of the cusp widths adds up to $w_1w_2$ (for nonprime $n$ 
some values of the twist parameters could correspond to nonprimitive 
surfaces, but for prime $n$ all surfaces with $n$ tiles are primitive).
All surfaces with $h_1$ and $h_2$ odd are in orbit $A$, all surfaces
with  $w_1$ and $w_2$ odd are in orbit $B$, and those with 
mixed parities for $w_i$ and $h_i$ are half in orbit $A$ half in 
orbit $B$.

\section{Asymptotics for a large prime number of squares}
\label{sec:asymptotics}

We need to estimate quantities of the type
\begin{equation*}
\tilde c(D_n)
= \frac{n}{\#D_n}
\sum
\frac{ab}{a^2}
\end{equation*}
where the sum is over positive integers $a$, $b$, $h$, $y$ satisfying
 conditions as above.

\subsection{A simpler sum}
\label{sec:simpler:sum}

Since $\#D$, for prime $n$, is asymptotically $\punfrac{3}{16}n^3$,
we first replace $\unfrac{n}{\#D_n}$ by $\unfrac{1}{n^2}$.
Second, we momentarily drop the parity conditions; we will reintroduce
them in the following subsections.
Last, we drop the condition $a < b$; we will explain later why this
does not change the asymptotic.

So we first consider the following simplified sum:
$$
S(n) = \sum_{a\geqslant 1} \frac{1}{a^2}
\sum_{b \geqslant 1}
\sum_{\substack{h \geqslant 1,\;y \geqslant 1 \\ a h + b y = n}}\frac{ab}{n^2}.
$$
Denote the sum over $b$ by $S(n,a)$.
Introducing the variable $m = by$,
$$
S(n,a) = \sum_{\substack{1\leqslant m \leqslant n - a \\ m \equiv n\!\pmod{a}}}
\sum_{b \mid m} \frac{ab}{n^2} = \frac{a}{n^2}
\cdot F(n - a, n, a)$$
$$
F(x,k,q) = \sum_{\substack{1 \leqslant m \leqslant x \\ m \equiv k
\pmod{q}}} \sum_{b\mid m}b. \leqno{\hbox{where}}
$$
The following asymptotics hold for $F(x,k,q)$, $S(n,a)$ and $S(n)$.

\begin{Lemma}
\label{lem:Fxkq}
For $k \wedge q = 1$ and $x \to \infty$,
$$
F(x,k,q) = \frac{x^2}{q} \cdot \frac{\pi^2}{12} \prod_{p\mid q}\Big(1 -
\frac{1}{p^2}\Big) + O_q(x \log x).
$$
\end{Lemma}
\begin{Lemma}
\label{lem:Sna} For $n$ prime, 
$\displaystyle
S(n,a) \rightarrow{} 
\frac{\pi^2}{12}\prod_{p\mid a} \Big(1 - \frac{1}{p^2}\Big)
$ as $n \rightarrow \infty$.
\end{Lemma}

\begin{Lemma}
\label{lem:Sn} For $n$ prime,
$\displaystyle
S(n) \rightarrow \frac{5}{4}$ as $n \to \infty$.
\end{Lemma}

\begin{proof}[Proof of \fullref{lem:Fxkq}]
If $m$ is prime to $q$, denote by $\bar m$ the integer in $\{0,
\ldots, q{-}1\}$ such that $\bar m m \equiv 1 \pmod{q}$, and by
$u = u(m,k,q)$ the integer in $\{0, \ldots, q{-}1\}$ such that $u
\equiv \bar m k \pmod{q}$; error terms depend on $q$.
\begin{equation*}
\begin{split}
F(x,k,q) 
& = \sum_{\substack{1 \leqslant md \leqslant x \\ md \equiv k \pmod{q}}} d
\\
&= \sum_{\substack{1 \leqslant m \leqslant x \\ m \wedge q = 1}}
\sum_{\substack{1 \leqslant d \leqslant x / m \\ d \equiv \bar m k
\pmod{q}}} d
\\
&= \sum_{\substack{1 \leqslant m \leqslant x \\ m \wedge q = 1}}
\sum_{\substack{1 \leqslant d \leqslant x / m \\ d \equiv u \pmod{q}}} d
\\
& = \sum_{\substack{1 \leqslant m \leqslant x \\ m \wedge q = 1}}
\sum_{1 \leqslant u + \lambda q \leqslant x / m} (u + \lambda q)
\\
&= \sum_{\substack{1 \leqslant m \leqslant x \\ m \wedge q = 1}}
\biggl(\Bigl(\sum_{1 \leqslant \lambda \leqslant \frac{1}{q}(\frac{x}{m} - u)}
\lambda q\Bigr) + O\Big(\frac{x}{m}\Big)\biggr)
\\
& = \sum_{\substack{1 \leqslant m \leqslant x \\ m \wedge q = 1}}
\Big(\frac{1}{2}q\Big(\frac{x}{qm}\Big)^2 + O\Big(\frac{x}{m}\Big) + O(1)\Big)
\\
&= \frac{x^2}{2q} 
\sum_{\substack{1 \leqslant m \leqslant x \\ m \wedge q = 1}}
\frac{1}{m^2} + O (x \log x)
\end{split}
\end{equation*}
To sum only over the integers $m$ with $m\wedge q = 1$, we
can sum over all $m$ with a factor $\mu(m\wedge q)$, so that all terms
cancel out except the ones we want.
$$\eqalignbot{
F(x,k,q) 
& = \frac{x^2}{2q} 
\sum_{d \mid q}\Big(\frac{\mu(d)}{d^2}\sum_{m\leqslant x /
d}\frac{1}{m^2}\Big) + O(x \log x)
\cr
&= \frac{x^2}{2q} 
\sum_{d \mid q}\frac{\mu(d)}{d^2}\Big(\frac{\pi^2}{6} + O(1/x)\Big) + O (x \log x)
\cr
& = \frac{x^2}{q} \cdot \frac{\pi^2}{12} 
\prod_{p\mid q}\Big(1 - \frac{1}{p^2}\Big) + O(x \log x)
}
\proved
$$
\end{proof}

\begin{proof}[Proof of \fullref{lem:Sna}]
The limit follows immediately from \fullref{lem:Fxkq} by a
dominated convergence argument (similar arguments were used in
\cite[\S\,7]{HL1}).
\end{proof}

\begin{proof}[Proof of \fullref{lem:Sn}]
This is a consequence of \fullref{lem:Sna} by the
following observation.
$$
\eqalignbot{
\sum_{a\geqslant 1}\frac{1}{a^2}\prod_{p \mid a}\Big(1 - \frac{1}{p^2}\Big)
& = \prod_{p}(1 + \sum_{\nu\geqslant 1}p^{-2\nu}(1 - p^{-2\nu}))
= \prod_{p}(1+p^{-2})
\cr
& = \prod_p\frac{1- p^{-4}}{1 - p^{-2}}
= \frac{\zeta(2)}{\zeta(4)}
= \frac{\pi^2/6}{\pi^4/90} = \frac{15}{\pi^2}
}\proved
$$
\end{proof}

\subsection{Sums with specified parities}

We introduce subsums of $S(n)$ for specified parities of the
parameters.

The observation we just made will need to be completed by the
following one.
$$
\sum_{\substack{a\geqslant 1\\ a\text{ even}}}
\frac{1}{a^2}\prod_{p \mid a}\Big(1 - \frac{1}{p^2}\Big)
= \sum_{\substack{a\geqslant 1\\ a\text{ odd}}}
\frac{1}{4a^2}\frac{3}{4}\prod_{p \mid a}\Big(1 - \frac{1}{p^2}\Big)
+ \sum_{\substack{a\geqslant 1\\ a\text{ even}}}
\frac{1}{4a^2}\prod_{p \mid a}\Big(1 - \frac{1}{p^2}\Big)
$$
$$
\sum_{\substack{a\geqslant 1\\ a\text{ odd}}}
\frac{1}{a^2}\prod_{p \mid a}\Big(1 - \frac{1}{p^2}\Big)
= \frac{12}{\pi^2} \quad
\text{and} \quad
\sum_{\substack{a\geqslant 1\\ a\text{ even}}}
\frac{1}{a^2}\prod_{p \mid a}\Big(1 - \frac{1}{p^2}\Big)
= \frac{3}{\pi^2}
\leqno{\hbox{so that}}
$$

\subsubsection{Odd widths}

We now consider the sum over odd $a$ and $b$:
$$
S^\text{ow}(n)
= \sum_{\substack{a\geqslant 1 \\ a\text{ odd}}} \frac{1}{a^2}
\sum_{\substack{b \geqslant 1 \\ b\text{ odd}}}
\sum_{\substack{h \geqslant 1, y \geqslant 1 \\ a h + b y = n}}
\frac{ab}{n^2}
$$
We proceed as for the sum $S(n)$, putting
$$\eqalign{
& F^\text{ow}(x,k,q)
= \sum_{\substack{1 \leqslant m \leqslant x \\ m \equiv k
\pmod{q}}} \sum_{\substack{b\mid m \\ b\text{ odd}}}b, \cr
& S^\text{ow}(n,a)
= \frac{a}{n^2} \cdot F^\text{ow}(n - a, n, a)
\quad \text{and} \quad
S^\text{ow}(n) 
= \sum_{\substack{a\geqslant 1\\a\text{ odd}}}\frac{1}{a^2} 
S^\text{ow}(n, a).
}$$

\begin{Lemma}
The following asymptotics hold for $F^\text{ow}(x,k,q)$,
$S^\text{ow}(n,a)$ and $S^\text{ow}(n)$.
\begin{gather*}
\tag*{\hbox{For odd $q$, odd $k$, and $x \to \infty$,\quad}}
F^\text{ow}(x,k,q) 
= \frac{x^2}{q}\frac{\pi^2}{24}\prod_{p\mid q}\Big(1 - \frac{1}{p^2}\Big)
+ O(x\log x).\\
S^\text{ow}(n,a)
\xrightarrow[\substack{n \to \infty\\n\text{ prime}}]{}
\frac{\pi^2}{24}\prod_{p \mid a}\Big(1 - \frac{1}{p^2}\Big).
\tag*{\hbox{For odd $a$,}}\\
S^\text{ow}(n)
\xrightarrow[\substack{n \to \infty\\n\text{ prime}}]{}
\frac{1}{2}.
\tag*{\hbox{Finally,}}
\end{gather*} 
\end{Lemma}

\begin{proof}
$$\eqalignbot{
F^\text{ow}(x,k,q) 
& = \sum_{t\geqslant 0} \sum_{\substack{1\leqslant m \leqslant x / 2^t
\cr 2^t m \equiv k \pmod q \\ m \equiv 1 \pmod 2}} \sum_{b\mid m}b
\cr
& = \sum_{t\geqslant 0} \biggl(\frac{(x/2^t)^2}{2q}\frac{\pi^2}{12}
\prod_{p \mid 2q}\Big(1 - \frac{1}{p^2}\Big)+O((x/2^t)\log(x/2^t))\biggr)
\cr
& = \frac{x^2}{q} \frac{1}{1 - \unfrac{1}{4}}\frac{\pi^2}{24}
\Big(1 - \frac{1}{2^2}\Big)\prod_{p \mid q}\Big(1 - \frac{1}{p^2}\Big)+O(x\log x)
\cr
& = \frac{x^2}{q}\frac{\pi^2}{24}\prod_{p\mid q}\Big(1 - \frac{1}{p^2}\Big)
+ O(x\log x)
}\proved
$$
\end{proof}

\subsubsection{Odd heights}

We now consider the sum over odd $h$ and $y$:
$$
S^\text{oh}(n)
= \sum_{a\geqslant 1} \frac{1}{a^2}
\sum_{b \geqslant 1}
\sum_{\substack{h \geqslant 1, y \geqslant 1 \\
h\text{, }y\text{ odd}
\\ a h + b y = n}}
\frac{ab}{n^2}.
$$
Proceeding as previously, we are led to introduce
$$
F^\text{oh}(x,k,q)
= \sum_{\substack{1 \leqslant m \leqslant x \\ m \equiv k + q
\pmod{2q}}} \sum_{\substack{b\mid m \\ m/b\text{ odd}}}b
\quad \text{and} \quad S^\text{oh}(n,a)
= \frac{a}{n^2} \cdot F^\text{oh}(n - a, n, a),
$$
$$ S^\text{oh}(n) 
= \sum_{a\geqslant 1}\frac{1}{a^2} S^\text{oh}(n, a).
\leqno{\hbox{and to write}}
$$

\begin{Lemma} 
The following asymptotics hold for $F^\text{oh}(x,k,q)$,
$S^\text{oh}(n,a)$ and $S^\text{oh}(n)$.
\begin{gather*}
\tag*{\hbox{For even $q$, odd $k$, and $x \to \infty$,\quad}}
F^\text{oh}(x,k,q) 
= \frac{x^2}{q}\frac{\pi^2}{24}\prod_{p\mid q}\Big(1 - \frac{1}{p^2}\Big)
+ O(x\log x).
\\
\tag*{\hbox{For odd $q$, odd $k$, and $x \to \infty$,\quad}}
F^\text{oh}(x,k,q) 
= \frac{x^2}{q}\frac{\pi^2}{32}\prod_{p\mid q}\Big(1 - \frac{1}{p^2}\Big)
+ O(x\log x).\\
S^\text{oh}(n,a)
\xrightarrow[\substack{n \to \infty\\n\text{ prime}}]{}
\frac{\pi^2}{24}\prod_{p \mid a}\Big(1 - \frac{1}{p^2}\Big).
\tag*{\hbox{For even $a$,}}\\
S^\text{oh}(n,a)
\xrightarrow[\substack{n \to \infty\\n\text{ prime}}]{}
\frac{\pi^2}{32}\prod_{p \mid a}\Big(1 - \frac{1}{p^2}\Big).
\tag*{\hbox{For odd $a$,}}\\
S^\text{oh}(n) \xrightarrow[\substack{n \to \infty\\n\text{ prime}}]{}
\frac{1}{2}.
\tag*{\hbox{Finally,}}
\end{gather*}
\end{Lemma}

\begin{proof}
For even $q$ and odd $k$:
\begin{equation*}
\begin{split}
F^\text{oh}(x,k,q) 
= \sum_{\substack{1 \leqslant m \leqslant x \\ m \equiv k + q
\pmod{2q}}} \sum_{b\mid m}b
& = \frac{x^2}{2q} \frac{\pi^2}{12}
\prod_{p \mid 2q}\Big(1 - \frac{1}{p^2}\Big) + O(x\log x)
\\
& = \frac{x^2}{q}\frac{\pi^2}{24}\prod_{p\mid q}\Big(1 - \frac{1}{p^2}\Big)
+ O(x\log x)
\end{split}
\end{equation*}
For odd $q$ and odd $k$:
$$\eqalignbot{
F^\text{oh}(x,k,q) 
& = \sum_{t\geqslant 1} \sum_{\substack{1\leqslant m \leqslant x / 2^t
\\ 2^t m \equiv k + q \pmod {2q} \\ m \text{ odd}}} \sum_{b\mid m} 2^t b
\cr
& = \sum_{t\geqslant 1} 2^t \sum_{\substack{1\leqslant m \leqslant x / 2^t
\cr 2^{t-1} m \equiv \upnfrac{k + q}{2} \pmod q \\ m \text{ odd}}} 
\sum_{b\mid m} b
\cr
& = \sum_{t\geqslant 1} 2^t
\frac{(x/2^t)^2}{2q}\frac{\pi^2}{12}\prod_{p\mid 2q}\Big(1 - \frac{1}{p^2}\Big)
+ O(x\log x)
\cr
& = \sum_{t \geqslant 1}\frac{1}{2^t} 
\frac{x^2}{q}\frac{\pi^2}{24}\Big(1 - \frac{1}{2^2}\Big )
\prod_{p\mid q}\Big(1 - \frac{1}{p^2}\Big) + O(x\log x)
\cr
& = \frac{x^2}{q}\frac{\pi^2}{32}
\prod_{p\mid q}\Big(1 - \frac{1}{p^2}\Big) + O(x\log x)
}\proved
$$
\end{proof}

\subsubsection{Mixed parities}

Dealing with the even--odd sums as above would be most cumbersome; 
this is fortunately not necessary.
Indeed,
$S(n) = S^\text{ow}(n) + S^\text{oh}(n) + S^\text{eo}(n)$,
and we know the limits of $S(n)$, $S^\text{ow}(n)$ and
$S^\text{oh}(n)$ when $n$ tends to infinity staying prime, so we have
$$
S^\text{eo}(n) 
\xrightarrow[\substack{n\to \infty \\n \text{ prime}}]{}
\frac{1}{4}.
$$

\subsection{Asymptotics for orbits A and B}

We end by showing that the obtained limit  is unchanged by
adding the condition $a<b$.

Indeed, since
$\#\{(h,y)\colon h\geqslant 1,\ y \geqslant 1,\ ah + b y = n\}
\leqslant n$, the sum
$$
\sum_{b = 1}^{a}
\sum_{\substack{h\geqslant 1,\;y \geqslant 1 \\ a h + b y = n}}
\frac{ab}{n^2}$$ 
is $O(1/n)$, where the constant of the $O$ depends on $a$.

Putting the previous sections together, $\tilde c(A_n)$ and $\tilde c(B_n)$
have the same asymptotics: $S^\text{A}(n) =
\punfrac{16}{3}(S^\text{oh}(n) + \frac{1}{2}S^\text{eo}(n))$ and
$S^\text{B}(n) = \punfrac{16}{3}(S^\text{ow}(n) +
\frac{1}{2}S^\text{eo}(n))$, so they both tend to $\unfrac{10}{3}$ as
$n$ tends to infinity, $n$ prime.


\section{Concluding remarks}

Numerical evidence suggests that the convergence to the generic
constant of the stratum occurs not only for prime $n$ but for general
$n$; however a proof would involve some complications in the
calculations which would make the exposition tedious.

A similar study for the constants that appear in the quadratic
asymptotics for the countings of saddle connections could also be
made.  There, one has to take into consideration both one-cylinder and
two-cylinder cusps, and some interesting phenomena can be observed.
Numerical calculations suggest that the sum of the contributions of
one-cylinder and two-cylinder cusps has a limit, but separate
countings for one-cylinder cusps do not have a limit for general $n$;
their asymptotics have fluctuations involving the prime factors of
$n$.

\bibliographystyle{gtart}
\bibliography{link}

\begin{thebibliography}{}
\providecommand\bibmarginpar{\leavevmode\marginpar}
\def\urlstyle#1{{\tt #1}}

\bibitem{DP}
\textbf{F Dal'bo}, \textbf{M Peign{\'e}}, \emph{Comportement asymptotique du
  nombre de g\'eod\'esiques ferm\'ees sur la surface modulaire en courbure non
  constante}, Ast\'erisque  (1996) 111--177 \xox{MR}{1634272}

\bibitem{EMWM}
\textbf{A Eskin}, \textbf{J Marklof}, \textbf{D Witte~Morris},
  \href{http://dx.doi.org/10.1017/S0143385705000234} {\emph{Unipotent flows on
  the space of branched covers of {V}eech surfaces}}, Ergodic Theory Dynam.
  Systems 26 (2006) 129--162 \xox{MR}{2201941}

\bibitem{EM}
\textbf{A Eskin}, \textbf{H Masur},
  \href{http://dx.doi.org/10.1017/S0143385701001225} {\emph{Asymptotic formulas
  on flat surfaces}}, Ergodic Theory Dynam. Systems 21 (2001) 443--478
  \xox{MR}{1827113}

\bibitem{EMS}
\textbf{A Eskin}, \textbf{H Masur}, \textbf{M Schmoll},
  \href{http://dx.doi.org/10.1215/S0012-7094-03-11832-3} {\emph{Billiards in
  rectangles with barriers}}, Duke Math. J. 118 (2003) 427--463
  \xox{MR}{1983037}

\bibitem{EMZ}
\textbf{A Eskin}, \textbf{H Masur}, \textbf{A Zorich},
  \href{http://dx.doi.org/10.1007/s10240-003-0015-1} {\emph{Moduli spaces of
  abelian differentials: the principal boundary, counting problems, and the
  {S}iegel--{V}eech constants}}, Publ. Math. Inst. Hautes \'Etudes Sci.  (2003)
  61--179 \xox{MR}{2010740}

\bibitem{F}
\textbf{C Faivre}, \emph{Distribution of {L}\'evy constants for quadratic
  numbers}, Acta Arith. 61 (1992) 13--34 \xox{MR}{1153919}

\bibitem{GuJu}
\textbf{E Gutkin}, \textbf{C Judge},
  \href{http://dx.doi.org/10.1215/S0012-7094-00-10321-3} {\emph{Affine mappings
  of translation surfaces: geometry and arithmetic}}, Duke Math. J. 103 (2000)
  191--213 \xox{MR}{1760625}

\bibitem{HL1}
\textbf{P Hubert}, \textbf{S Leli{\`e}vre},
  \href{http://dx.doi.org/10.1007/BF02777365} {\emph{Prime arithmetic
  {T}eichm\"uller discs in $\mathcal{H}(2)$}}, Israel J. Math. 151 (2006)
  281--321 \xox{MR}{2214127}

\bibitem{LeRo}
\textbf{S Leli\`evre}, \textbf{E Royer}, \emph{Orbitwise countings in
  $\mathcal{H}(2)$ and quasimodular forms}, to appear in Internat.\ Math.\
  Res.\ Notices  (2006)

\bibitem{Ma88}
\textbf{H Masur}, \emph{Lower bounds for the number of saddle connections and
  closed trajectories of a quadratic differential}, from: ``Holomorphic
  functions and moduli, Vol.\ I (Berkeley, CA, 1986)'', (D Drasin, editor),
  Math. Sci. Res. Inst. Publ. 10, Springer, New York (1988)  215--228
  \xox{MR}{955824}

\bibitem{Ma90}
\textbf{H Masur}, \emph{The growth rate of trajectories of a quadratic
  differential}, Ergodic Theory Dynam. Systems 10 (1990) 151--176
  \xox{MR}{1053805}

\bibitem{McDS}
\textbf{C\,T McMullen}, \href{http://dx.doi.org/10.1007/s00208-005-0666-y}
  {\emph{Teichm\"uller curves in genus two: discriminant and spin}}, Math. Ann.
  333 (2005) 87--130 \xox{MR}{2169830}

\bibitem{Schmo}
\textbf{M Schmoll}, \href{http://dx.doi.org/10.1007/s00039-002-8260-x}
  {\emph{On the asymptotic quadratic growth rate of saddle connections and
  periodic orbits on marked flat tori}}, Geom. Funct. Anal. 12 (2002) 622--649
  \xox{MR}{1924375}

\bibitem{Ve89}
\textbf{W\,A Veech}, \href{http://dx.doi.org/10.1007/BF01388890}
  {\emph{Teichm\"uller curves in moduli space, {E}isenstein series and an
  application to triangular billiards}}, Invent. Math. 97 (1989) 553--583
  \xox{MR}{1005006}

\bibitem{Vo97}
\textbf{Y\,B Vorobets},
  \href{http://dx.doi.org/10.1070/SM1997v188n03ABEH000211} {\emph{Ergodicity of
  billiards in polygons}}, Mat. Sb. 188 (1997) 65--112 \xox{MR}{1462024}

\end{thebibliography}

\end{document}